\documentclass{article}
\usepackage{amsfonts}
\usepackage{amssymb}

\newenvironment{proof}[1][Proof]{\noindent\textbf{#1.} }{\ \rule{0.5em}{0.5em}}
\input{tcilatex}

\begin{document}

\begin{center}
\begin{eqnarray*}
&&\text{{\LARGE \ }{\large A characterization of the distances between Pisot
numbers }} \\
&&\text{{\large \ \ \ \ \ \ \ \ \ \ \ \ \ \ \ \ \ \ \ \ generating the same
number field }}
\end{eqnarray*}

\smallskip

by M. J. BERTIN\ \ and \ T. ZA\"{I}MI

\bigskip

Ce papier est dedi\'{e} \`{a} la m\'{e}moire du Professeur Georges RHIN

\textit{\ }
\end{center}

\textbf{Abstract. \ }\textit{Let }$K$\textit{\ be a real algebraic number
field\ and let }$\wp _{K}$\textit{\ be the set of Pisot numbers generating }$%
K.$\textit{\ We show that the elements of }$\wp _{K}-\wp _{K}$\textit{\ are
the algebraic integers of }$K$\textit{\ whose images under the action of all
embeddings of }$K$\textit{\ into }$\mathbb{C},$\textit{\ other than the
identity of }$K,$\textit{\ are of modulus less than }$2.$\textit{\ This
completes certain previous results due to Dubickas and the second author.
Also, from the proof of this result, based on Minkowski's theorem on
compact, convex and symmetric bodies, we deduce some related known results.}

\bigskip


\medskip


\medskip

\textbf{1. Introduction}

\textbf{\ } Let throughout $K$ designate a real algebraic number field
having $s$ real conjugates and $2t$ non-real conjugates, where $s\geq 1,$ $%
t\geq 0$ and $d:=s+2t$ is the degree of $K.$ Let $\sigma _{1},...,\sigma
_{d} $ denote the distinct embeddings of $K$ into $\mathbb{C},$ labelled so
that $\sigma _{1}$ is the identity of $K,$ $\sigma _{j}$ is real if and only
if $j\leq s,$ and $\sigma _{s+t+j},$ where $j\in \{1,...,t\},$ is equal to
the composition of the complex conjugation and the embedding $\sigma _{s+j}.$

A Pisot number is a real algebraic integer greater than $1$ whose other
conjugates are of modulus less than $1,$ and the set of such numbers is
traditionally denoted by $S.$ Let 
\[
\wp _{K}:=\{\theta \mid \theta \in S\cap K\text{ \ and \ }K=\mathbb{Q}%
(\theta )\}. 
\]%
Pisot [9] was the first to investigate the set $\wp _{K}$ and he showed, in
particular, that $\wp _{K}$ contains units, whenever $d\geq 2.$ Since there
are at most a finite number of algebraic integers of fixed degree and having
all conjugates in a bounded subset of the complex plane, the set $\wp _{K}$
is a discrete subset of $\mathbb{R},$ i.e., any finite real interval
contains at most a finite number of such numbers. Moreover, if $\theta
^{\prime }>\theta $ \ are any elements of $\wp _{K},$ then 
\[
1\leq \left\vert Norm(\theta ^{\prime }-\theta )\right\vert =(\theta
^{\prime }-\theta )\dprod_{2\leq j\leq d}\left\vert \sigma _{j}(\theta
^{\prime })-\sigma _{j}(\theta )\right\vert \leq 2^{d-1}(\theta ^{\prime
}-\theta ), 
\]%
$\theta ^{\prime }-\theta \geq 2^{1-d},$ and so the set $\wp _{K}$ is
uniformly discrete; recall that a real set is said to be uniformly discrete
(or more precisely, $\rho $-discrete) if there is a positive number $\rho $
such that each real interval of the form $(r,r+\rho ]$ contains at most one
element of this set.

A less trivial property of the set $\wp _{K}$ is that it is also relatively
dense in $[1,\infty ),$ that is, there is a real number $\rho ^{\prime }>0$
such that every subinterval of $[1,\infty )$ of the form $(\varepsilon
,\varepsilon +\rho ^{\prime }]$ contains at least one element of $\wp _{K};$
to be more precise, we say in this case that $\wp _{K}$ is $\rho ^{\prime }$%
-dense in $[1,\infty ).$ This property is an immediate corollary of [4,
Theorem 1.4]. Using some results, due to Meyer, on harmonious sets (see for
instance [6], [7] and [8]), the second author pointed out, in [12], that the
family $\wp _{K}^{\prime }:=\wp _{K}\cup (-\wp _{K})$ is a real Meyer set,
i.e., there is a finite subset $F$ of the real line $\mathbb{R}$ such that $%
\wp _{K}^{\prime }-\wp _{K}^{\prime }\subset \wp _{K}^{\prime }+F,$ and a
generalization of this proposition for complex Pisot numbers was given in
[1]. A real Meyer set is, in particular, a relatively dense and uniformly
discrete subset of $\mathbb{R}.$

Let $\theta _{1}<\theta _{2}<\theta _{3}<\cdot \cdot \cdot $ \ denote the
elements of $\wp _{K}$ and 
\[
D_{K}:=\{\theta _{n}-\theta _{m}\mid (n,m)\in \mathbb{N}^{2}\text{ \ and }%
n>m\}. 
\]%
As signaled in [13], $(\wp _{K}-\wp _{K})\cap (0,\infty )=D_{K},$ $\wp
_{K}-\wp _{K}=(-D_{K})\cup \{0\}\cup D_{K}$ and 
\begin{equation}
D_{K}\subset E_{K}:=\{\beta \in \mathbb{Z}_{K}\cap (0,\infty )\mid
\left\vert \sigma _{2}(\beta )\right\vert <2,...,\left\vert \sigma
_{d}(\beta )\right\vert <2\},  \tag{1}
\end{equation}%
where $\mathbb{Z}_{K}$ denotes the ring of integers of $K.$

Also, [13, Theorem 1.1], [3, Theorem 2] and [13, Theorem 1.4] yield,
respectively,%
\[
1\in D_{K},\text{ \ \ }\{\beta \in E_{K}\mid K=\mathbb{Q}(\beta )\}\subset
D_{K}, 
\]%
and $E_{K}\subset D_{K}$ when\textit{\ }$K$\ is totally real.

Clearly, $\wp _{\mathbb{Q}}=\mathbb{N\diagdown \{}1\},$ $\wp _{\mathbb{Q}%
}-\wp _{\mathbb{Q}}=\mathbb{Z}$ \ and $D_{\mathbb{Q}}=\mathbb{N}=E_{\mathbb{Q%
}}.$ The aim of this note is to show, independently of the above mentioned
results, that $E_{K}\subset D_{K}$ for any real number field\textit{\ }$K.$
This assertion is a consequence of the following theorem.

\medskip

\textbf{Theorem 1. }\textit{\ With the notation above, where }$K$\textit{\
denotes a real number field of degree }$s+2t\geq 2,$ \textit{let }$%
(x_{2},...,x_{s+t})$\textit{\ be an element of the Cartesian product} $%
\mathbb{R}^{s-1}\mathbb{\times C}^{t}$ \textit{whose all components are of
modulus less than} $1,$ \textit{and} \textit{let }$\varepsilon $\textit{\ be
a sufficiently small positive number such that all real intervals }$%
[x_{2}-\varepsilon ,x_{2}+\varepsilon ],...,$\textit{\ }$[x_{s}-\varepsilon
,x_{s}+\varepsilon ]$\textit{\ and all complex discs }$\{z\in \mathbb{C}\mid
\left\vert z-x_{s+1}\right\vert \leq \varepsilon \},...,$\textit{\ }$\{z\in 
\mathbb{C}\mid \left\vert z-x_{s+t}\right\vert \leq \varepsilon \}$\textit{\
are contained in the open unit disc} $\{z\in \mathbb{C}\mid \left\vert
z\right\vert <1\}.$

\textit{Then, there is a positive constant }$c=c(K,\varepsilon )$ \textit{%
such that any real interval }$[x_{1}-c,x_{1}+c]$ \textit{contains an
algebraic integer }$\theta =\theta (K,\varepsilon ,x_{1})\in K$\textit{\ \
such that }$\left\vert x_{j}-\sigma _{j}(\theta )\right\vert \leq
\varepsilon $ \textit{for any} $j\in \{2,...,s+t\}.$

\bigskip

\textbf{Corollary 1.}\textit{\ The equality }$E_{K}=D_{K}$\textit{\ holds
for any real number field} $K.$

\bigskip

As mentioned above, the inclusion $E_{K}\subset D_{K},$ where $K$ is totally
real, has been shown in [13, Theorem 1.4] and its proof is based on [13,
Lemma 3.4], saying that the set $\{(\sigma _{2}(\theta ),...,\sigma
_{d}(\theta ))$ $\mid \theta \in \wp _{K}\}$ is dense\ in the hypercube $%
[-1,1]^{d-1}.$ From Theorem 1 we can also deduce a generalization of this
last mentioned lemma.

\medskip

\textbf{Corollary 2. }\textit{\ Let }$K$\textit{\ be a real number field
with degree greater than }$1.$ \textit{Then, the set }$\{(\sigma _{2}(\theta
),...,\sigma _{s+t}(\theta ))$ $\mid \theta \in \wp _{K}\}$\textit{\ is
dense in} $[-1,1]^{s-1}\times \{z\in \mathbb{C}\mid \left\vert z\right\vert
\leq 1\}^{t}.$

\bigskip

Let $C_{K}:=\{\theta _{n+1}-\theta _{n}\mid \mathit{\ }n\in \mathbb{N}\}$ be
the set of positive differences of consecutive elements of $\wp _{K}.$ Then,
[13, Theorem 1.7] states that $C_{K}$ is finite, $\func{card}(C_{K})\geq 2,$
and\textit{\ }$\func{card}(C_{K})\geq 2^{d-1}$ whenever $K$\ is totally
real. In fact, Corollary 2 implies the following more general result.

\medskip \medskip

\textbf{Corollary 3.} \textit{We have} $\func{card}(C_{K})\geq 2^{s+t-1}$ 
\textit{for any number field} $K.$

\medskip

Let $\epsilon \in (0,1].$ Following [4], a real algebraic integer $\theta >1$
whose other conjugates are of modulus less than $\epsilon $ is called an $%
\epsilon $-Pisot number (for $\epsilon =1$ we obtain the usual Pisot
numbers). Some properties of $\epsilon $-Pisot numbers may be found in [11].
Also, we can deduce from Theorem 1 another proof of the following result,
due to Fan and Schmeling.

\medskip

\textbf{Corollary 4. [4]}\textit{\ The set of }$\epsilon $-\textit{Pisot
numbers, generating the real number field }$K,$\textit{\ is relatively dense
in }$[1,\infty ).$

\bigskip

The proofs of all these results are given in the next section. The proof of
Theorem 1 is based on Minkowski's theorem, saying that if $E$ is a compact,
convex and symmetric subset of a $d$-dimensional $\mathbb{R}$-vector space $%
V,$ with volume $\limfunc{vol}(E)=2^{d}\limfunc{vol}(\Lambda ),$ where $%
\Lambda $ is a lattice in $V,$ then $E\cap \Lambda $ contains a non-zero
point (see for instance [5, pages 159-160] or [2, page 51]).

Throughout, when we speak about conjugates, the norm and the degree of an
algebraic number, without mentioning the basic field, this is meant over $%
\mathbb{Q}.$ Also, the degree, the discriminant and the conjugates of a
number field are considered over $\mathbb{Q},$ and an integer means a
rational integer.

\bigskip

\bigskip

\textbf{2. The proofs}

\smallskip

\textbf{Proof of Theorem 1.\ }Consider the map $\sigma $ from the field $K$
to the $d$-dimensional $\mathbb{R}$-vector space $V:=\mathbb{R}^{s}\mathbb{%
\times C}^{t},$ defined by 
\[
\sigma (\gamma )=(\sigma _{1}(\gamma ),...,\sigma _{s+t}(\gamma )), 
\]%
where $\gamma $ is any element of $K$ ($\sigma $ associates to $\gamma $ all
its real conjugates and exactly one from each complex conjugate pair of
non-real conjugates).

Then, $\sigma $ is an injective additive group homomorphism (in fact, $%
\sigma $ is also a ring homomorphism and the field $\sigma (K)$ is a subring
of $V),$ and the image $\Lambda :=$ $\sigma (\mathbb{Z}_{K})$ of the ring $%
\mathbb{Z}_{K}$ of the integers of $K$ is a subgroup of $V.$

Fix a basis $\{\omega _{1},...,\omega _{d}\}$ of the $\mathbb{Z}$-module $%
\mathbb{Z}_{K}.$ Then, the vectors 
\[
v_{1}:=\sigma (\omega _{1}),\text{ }...,\text{ }v_{d}:=\sigma (\omega _{d})
\]%
form a basis for the $\mathbb{Z}$-module $\Lambda .$ In order to show that $%
\Lambda $ is a lattice of $V,$ it suffices to show that the $d$ vectors $%
v_{1},...,v_{d}$ are $\mathbb{R}$-independent (and so generate the $\mathbb{R%
}$-vector space $V).$ To do this, consider the map $\tau $ from $V$ to the $d
$-dimensional $\mathbb{R}$-vector space $\mathbb{R}^{d},$ defined by $\tau
(x_{1},...,x_{s+t})=$  \ 
\[
(x_{1},...,x_{s},\frac{x_{s+1}+\overline{x_{s+1}}}{2},\frac{x_{s+1}-%
\overline{x_{s+1}}}{2i},...,\frac{x_{s+t}+\overline{x_{s+t}}}{2},\frac{%
x_{s+t}-\overline{x_{s+t}}}{2i}),
\]%
where $i^{2}=-1$ and the bar denotes the complex conjugation. Then, $\tau $
is an injective linear transformation, and so it is an isomorphism. Suppose
that $\mathbb{R}^{d}$ (resp. $V)$ is endowed with the canonical basis (resp.
with the usual topology), and let $\delta $ denote the determinant of the $d$
vectors $\tau (v_{j})=$ 
\[
(\sigma _{1}(\omega _{j}),...,\sigma _{s}(\omega _{j}),\func{Re}(\sigma
_{s+1}(\omega _{j})),\func{Im}(\sigma _{s+1}(\omega _{j})),...,\func{Re}%
(\sigma _{s+t}(\omega _{j})),\func{Im}(\sigma _{s+t}(\omega _{j}))),
\]%
where $j\in \{1,...,d\}.$ Then, simple linear combinations on the rows of $%
\delta ,$ corresponding to non-real embeddings of $K$ into $\mathbb{C},$
followed by permutations of these rows yield 
\begin{equation}
\delta =(-1)^{t(t-1)/2}\frac{\det [\sigma _{j}(\omega _{k})]_{1\leq j,k\leq
d}}{(-2i)^{t}}.  \tag{2}
\end{equation}%
Because $(\det [\sigma _{j}(\omega _{k})]_{1\leq j,k\leq d})^{2}$ is equal
to the discriminant, say $\Delta _{K},$ of the field $K,\ \ $we have, by
(2), that $\delta \neq 0,$ the vectors $\tau (v_{1}),...,\tau (v_{d})$ are $%
\mathbb{R}$-linearly independent, and so are the elements $v_{1},...,v_{d}$
of $V.$ Hence, $\Lambda $ is a lattice of $V,$ $\tau (\Lambda )$ is a
lattice of $\mathbb{R}^{d}$ with volume $\limfunc{vol}(\tau (\Lambda
))=\left\vert \delta \right\vert ,$ and so, by (2), 
\[
\limfunc{vol}(\Lambda )=\frac{\sqrt{\left\vert \Delta _{K}\right\vert }}{%
2^{t}}
\]%
(for more related details see [2, pages 156-157]).

Let $\rho _{1},...,\rho _{s+t}$ be sufficiently large positive numbers so
that the fundamental parallelepiped $\{\epsilon _{1}v_{1}+\cdot \cdot \cdot
+\epsilon _{d}v_{d}\mid (\epsilon _{1},...,\epsilon _{d})\in \lbrack
0,1]^{d}\}$ of $\Lambda $ is contained in the closed ball 
\[
\text{\ss }:=\dprod\limits_{1\leq j\leq s}[-\rho _{j},\rho _{j}]\times
\dprod\limits_{1\leq j\leq t}\{z\in \mathbb{C}\mid \left\vert z\right\vert
\leq \sqrt{\rho _{s+j}}\} 
\]%
of $V.$ Then, any element $v$ of $V$ may be written $v=r_{1}v_{1}+\cdot
\cdot \cdot +r_{d}v_{d}$ for some $(r_{1},...,r_{d})\in \mathbb{R}^{d},$ and
so 
\[
v=([r_{1}]v_{1}+\cdot \cdot \cdot +[r_{d}]v_{d})+((r_{1}-[r_{1}])v_{1}+\cdot
\cdot \cdot +(r_{d}-[r_{d}])v_{d}), 
\]%
where $[.]$ is the integer part function, $[r_{1}]v_{1}+\cdot \cdot \cdot
+[r_{d}]v_{d}\in \Lambda $ and $(r_{1}-[r_{1}])v_{1}+\cdot \cdot \cdot
+(r_{d}-[r_{d}])v_{d}\in $\ss $;$ thus 
\begin{equation}
V\subset \Lambda +\text{\ss }.  \tag{3}
\end{equation}%
Let $(x_{1},...,x_{s+t})$ be any element of $V$ such that $\left\vert
x_{j}\right\vert <1$ for any $j\in \{2,...,s+t\},$ and let $\varepsilon $ be
a sufficiently small positive number such that all real intervals $%
[x_{2}-\varepsilon ,x_{2}+\varepsilon ],...,$\textit{\ }$[x_{s}-\varepsilon
,x_{s}+\varepsilon ]$\textit{\ }and all discs\textit{\ }$\{z\in \mathbb{C}%
\mid \left\vert z-x_{s+1}\right\vert \leq \varepsilon \},...,$\textit{\ }$%
\{z\in \mathbb{C}\mid \left\vert z-x_{s+t}\right\vert \leq \varepsilon \}$%
\textit{\ }are contained in the open unit disc $\{z\in \mathbb{C}\mid
\left\vert z\right\vert <1\}.$

Set%
\begin{equation}
c:=\frac{2^{t}\sqrt{\left\vert \Delta _{K}\right\vert }\rho _{1}\cdot \cdot
\cdot \rho _{s+t}}{\pi ^{t}\varepsilon ^{d-1}},  \tag{4}
\end{equation}%
and

\[
E:=[-\frac{c}{\rho _{1}},\frac{c}{\rho _{1}}]\times \dprod\limits_{2\leq
j\leq s}[-\frac{\varepsilon }{\rho _{j}},\frac{\varepsilon }{\rho _{j}}%
]\times \dprod\limits_{1\leq j\leq t}\{z\in \mathbb{C}\mid \left\vert
z\right\vert \leq \frac{\varepsilon }{\sqrt{\rho _{s+j}}}\}.
\]%
Since $(\rho _{1},...,\rho _{s+t})$ depends only on the basis $\{\omega
_{1},...,\omega _{d}\}$ of $\mathbb{Z}_{K},$ and $\varepsilon $ does not
depend on $x_{1},$ we see that $c$ depends only on the pair $(K,\varepsilon
).$\ Also, $E$ is a compact, convex and symmetric (about the origin) subset
of $V$ \ with volume%
\[
\limfunc{vol}(E)=\frac{2^{s}\pi ^{t}\varepsilon ^{d-1}c}{\rho _{1}\cdot
\cdot \cdot \rho _{s+t}}=2^{s+t}\sqrt{\left\vert \Delta _{K}\right\vert }%
=2^{d}\limfunc{vol}(\Lambda ).
\]%
It follows, by Minkowski's theorem, that $\Lambda \cap E$ contains a
non-zero point, i.e., there is $\alpha \in \mathbb{Z}_{K}$ such that $\alpha
\neq 0,$ 
\[
\left\vert \alpha \right\vert \leq \frac{c}{\rho _{1}},\text{ \ }\left\vert
\sigma _{2}(\alpha )\right\vert \leq \frac{\varepsilon }{\rho _{2}},\text{ }%
...,\text{ }\left\vert \sigma _{s}(\alpha )\right\vert \leq \frac{%
\varepsilon }{\rho _{s}}
\]%
and 
\[
\left\vert \sigma _{s+1}(\alpha )\right\vert \leq \frac{\varepsilon }{\sqrt{%
\rho _{s+1}}},\text{ }...,\text{ }\left\vert \sigma _{s+t}(\alpha
)\right\vert \leq \frac{\varepsilon }{\sqrt{\rho _{s+1}}}.
\]%
Finally, applying (3) to the vector\ $(x_{1}/\alpha ,$ $x_{2}/\sigma
_{2}(\alpha ),...,$ $x_{s+t}/\sigma _{s+t}(\alpha ))$ of $V,$ we get \ 
\[
\frac{x_{1}}{\alpha }=\beta +b_{1},\text{ \ }\frac{x_{2}}{\sigma _{2}(\alpha
)}=\sigma _{2}(\beta )+b_{2},\text{ }...,\text{ }\frac{x_{s+t}}{\sigma
_{s+t}(\alpha )}=\sigma _{s+t}(\beta )+b_{s+t},
\]%
for some $\beta \in \mathbb{Z}_{K}$ and $b=(b_{1},...,b_{s+t})\in $\ss $,$
and so $\left\vert x_{1}-\theta \right\vert \leq c,$%
\[
\text{ }\left\vert x_{2}-\sigma _{2}(\theta )\right\vert =\left\vert
b_{2}\sigma _{2}(\alpha )\right\vert \leq \varepsilon ,\text{ }...,\text{ }%
\left\vert x_{s+t}-\sigma _{s+t}(\theta )\right\vert =\left\vert
b_{s+t}\sigma _{s+t}(\alpha )\right\vert \leq \varepsilon ,
\]%
where $\theta :=\alpha \beta \in $ $\mathbb{Z}_{K}.$%
\endproof%

\medskip \medskip

\textbf{Proof of Corollary 1. }\ As noted above, we have $\ E_{\mathbb{Q}}=$ 
$D_{\mathbb{Q}}.$ Suppose $d\geq 2.$ Then, the relation (1) says that $%
D_{K}\subset $ $E_{K}.$ To show the inverse inclusion $E_{K}\subset $ $%
D_{K}, $ consider an element $\beta $ of $E_{K}.$ Then, $\beta \in \mathbb{Z}%
_{K}\cap (0,\infty )$ and 
\[
\rho :=\max \{\left\vert \sigma _{2}(\beta )\right\vert ,...,\left\vert
\sigma _{s+t}(\beta )\right\vert \}<2. 
\]%
Suppose $(x_{2},...,x_{s+t}):=(\sigma _{2}(\beta )/2,...,\sigma _{s+t}(\beta
)/2)$ and 
\[
\varepsilon \in (0,1-\rho /2) 
\]
satisfy the assumption of Theorem 1. Clearly, Theorem 1 gives that the
interval $[2\beta +1,2\beta +1+2c]$ contains an algebraic integer $\theta
=\theta (K,\varepsilon ,2\beta +1+c)\in K$\ such that $\left\vert \sigma
_{j}(\beta )-2\sigma _{j}(\theta )\right\vert \leq 2\varepsilon $ for any $%
j\in \{2,...,s+t\}.$

Because 
\[
\theta >\theta -\beta \geq 1+\beta >1 
\]%
and all the numbers $\sigma _{2}(\theta ),...,\sigma _{s+t}(\theta )$ belong
to the open unit disc (this follows from the hypothesis of Theorem 1), we
see that $\theta \in \wp _{K}.$ \ On the other hand, we have (by the same
way as in the proof of [3, Theorem 2]) that 
\[
2\left\vert \sigma _{j}(\beta )-\sigma _{j}(\theta )\right\vert \leq
\left\vert \sigma _{j}(\beta )-2\sigma _{j}(\theta )+\sigma _{j}(\beta
)\right\vert \leq 2\varepsilon +\rho <2, 
\]%
for any $j\in \{2,...,s+t\},$ and so $\left\vert \sigma _{j}(\theta -\beta
)\right\vert =\left\vert \sigma _{j}(\theta )-\sigma _{j}(\beta )\right\vert
<1.$ Hence, $\theta -\beta \in \wp _{K}$ and $\beta =\theta -(\theta -\beta
)\in D_{K}.$ Finally, notice if we consider the intervals $[(n-1)c,(n+1)c],$
where $n$ is sufficiently large integer, instead of the interval $[2\beta
+1,2\beta +1+2c],$ we obtain that any $\beta \in E_{K}$ can be written as a
difference of two elements of $\wp _{K}$ in infinite number of ways.%
\endproof%

\bigskip

\textbf{Proof of Corollary 2. }Let $(x_{2},...,x_{s+t})$ be a vector of $%
\mathbb{R}^{s-1}\mathbb{\times C}^{t}$ whose components are all of modulus
less than $1,$ and let $\varepsilon >0.$\textit{\ }To show the desired
result we may suppose, without loss of generality, that $\varepsilon $
satisfies the assumption of Theorem 1. Then, this theorem gives that for any
real number $x_{1}>c+1$ there is $\theta \in \mathbb{Z}_{K}\cap (1,\infty )$%
\textit{\ }such that $\left\vert x_{j}-\sigma _{j}(\theta )\right\vert \leq
\varepsilon $ for any $j\in \{2,...,s+t\}.$ Since all numbers $\sigma
_{2}(\theta ),...,\sigma _{s+t}(\theta )$ belong to the open unit disc, we
immediately deduce that $\theta \in \wp _{K}$ and Corollary 2 is true.%
\endproof%

\bigskip

\textbf{Proof of Corollary 3. \ }Since $C_{\mathbb{Q}}=\{1\},$ Corollary 3
is true for $d=1.$ Suppose $d\geq 2.$ Then, $s+t-1\geq 1,$ and we proceed by
the same way as in the proof of [13, Theorem 1.7(ii)].

Fix an element $(\epsilon _{2},...,\epsilon _{s+t})$ of the set $%
\{-1,1\}^{s+t-1}.$ Then, Corollary 2 gives that there is a subsequence $%
(\theta _{f(n)})_{n\in \mathbb{N}}$ of the sequence $(\theta _{n})_{n\in 
\mathbb{N}}$ such that%
\[
\lim_{n\rightarrow \infty }(\sigma _{2}(\theta _{f(n)}),...,\sigma
_{s+t}(\theta _{f(n)}))=(\epsilon _{2},...,\epsilon _{s+t}). 
\]%
Because $\{\theta _{f(n)+1}-\theta _{f(n)}\mid $ $n\in \mathbb{N\}}$ is
contained in the finite set $C_{K},$ there exist an element $c$ of $C_{K}$
and an infinite subset $N$ of $\mathbb{N}$ such that $\theta
_{f(n)+1}-\theta _{f(n)}=c$ \ for all $n\in N.$ Therefore, for each $%
(j,n)\in \{2,...,s+t\}\times N,$ $\ \sigma _{j}(\theta _{f(n)})+\sigma
_{j}(c)=\sigma _{j}(\theta _{f(n)+1}),$ and hence 
\[
\lim_{n\rightarrow \infty }(\sigma _{j}(\theta _{f(n)+1})=\epsilon
_{j}+\sigma _{j}(c)=(\epsilon _{j}+\func{Re}(\sigma _{j}(c)))+i\func{Im}%
(\sigma _{j}(c)). 
\]%
It follows that 
\begin{equation}
(\epsilon _{j}+\func{Re}(\sigma _{j}(c)))^{2}+\func{Im}(\sigma
_{j}(c))^{2}\leq 1,  \tag{5}
\end{equation}%
for any $j\in \{2,...,s+t\},$ $\func{Re}(\sigma _{j}(c))<0$ (resp. $\func{Re}%
(\sigma _{j}(c))>0)$ when $\varepsilon _{j}=1$ (resp. $\varepsilon _{j}=-1),$
since the equality $\func{Re}(\sigma _{j}(c))=0,$ gives, by (5), $\func{Im}%
(\sigma _{j}(c))=0$ and $c=0.$

Consequently, for any fixed $(\epsilon _{2},...,\epsilon _{s+t})\in
\{-1,1\}^{s+t-1}$ there is $c\in C_{K}$ such that for each $j\in
\{2,...,s+t\}$ the numbers $\func{Re}(\sigma _{j}(c))$ and $\varepsilon _{j}$
have opposite signs; thus there are $2^{s+t-1}$ distinct elements of the set 
$C_{K}$ which are in a one-to-one correspondence with the elements of $%
\{-1,1\}^{s+t-1}.$%
\endproof%

\bigskip

\textbf{Proof of Corollary 4. }Clearly, $\wp _{Q}$ is $1$-dense. Suppose $%
d\geq 2,$ $\epsilon \in (0,1],$ $\varepsilon \in (0,\epsilon )$ and $%
(x_{2},...,x_{s+t}):=(0,...,0).$ Then, $\varepsilon <1,$ and we have by
Theorem 1 that\textit{\ }any real interval $[r,r+2c]\subset (0,\infty )$
contains an element $\theta =\theta (K,\varepsilon ,r+c)$ of $\mathbb{Z}_{K}$
such that $\left\vert \sigma _{j}(\theta )\right\vert \leq \varepsilon $ for
any $j\in \{2,...,s+t\}.$ Since the norm of the positive algebraic $\theta $
is an integer, and the other conjugates of $\theta $ are of modulus less
than $1,$ we have $\theta >1$ and so $\theta $ is not repeated under the
action of the embeddings of $K$ into $\mathbb{C}.$ Hence, $\theta $ is an $%
\epsilon $-Pisot, satisfying $K=\mathbb{Q}(\theta ),$ and the set of $%
\epsilon $-Pisot numbers, generating $K,$ is a $2c$-dense subset of $%
(0,\infty ).$%
\endproof%

\bigskip

\textbf{Remark. }From the proof of Corollary 4 we deduce that $\min \wp
_{K}\leq 2c$ and $\max C_{K}\leq 2c,$ where $c$ is given by the relation
(4). In fact, by letting $\varepsilon $ tend to $1,$ we obtain that 
\[
B_{K}:=\frac{2^{t}\sqrt{\left\vert \Delta _{K}\right\vert }\rho _{1}\cdot
\cdot \cdot \rho _{s+t}}{\pi ^{t}} 
\]%
is a upper bound for $\min \wp _{K}$ and for $\max C_{K}.$ Some numerical
investigations suggest that this bound is far from being optimal. For
instance, if we consider a quadratic field $K=\mathbb{Q}(\sqrt{m}),$ where $%
m\equiv 2,3$ $\func{mod}4$ is a square-free integer, and we set $(\omega
_{1},\omega _{2},\rho _{1},\rho _{2})=(1,\sqrt{m},1+\sqrt{m},\sqrt{m})$ in
the proof of Theorem 1, then we get $B_{K}:=\sqrt{\left\vert \Delta
_{K}\right\vert }\rho _{1}\rho _{2}=4m(1+\sqrt{m}).$ On the other hand, the
computation performed in [13, Section 2] gives $\min \wp _{K}=\max
C_{K}=\left\lfloor \sqrt{m}\right\rfloor +\sqrt{m},$ where $\left\lfloor
.\right\rfloor $ is the integer part function. Finally, note that V\'{a}vra
and Veneziano dealt in [10] with the algorithmic problem of finding $\min
\wp _{K}$ and some related questions.

\bigskip \medskip\ 

\textbf{References}

[1] M. J. Bertin and T. Za\"{\i}mi, \textit{Complex Pisot numbers in
algebraic number fields}, C. R. Math., Acad. Sci. Paris \textbf{353}, No. 11
(2015), 965-967.

[2] R. Descombes, \textit{\'{E}l\'{e}ments de th\'{e}orie des nombres},
Presses Universitaires de France, 1986.

[3] A. Dubickas, \textit{Numbers expressible as a difference of two Pisot
numbers}, Acta Math. Hung. \textbf{172} (2024), 346-358.

[4] A. H. Fan and J. Schmeling, $\varepsilon -$\textit{Pisot numbers in any
real algebraic number field are relatively dense}, J. Algebra \textbf{272}
(2004), 470-475.

[5] A. Fr\"{o}hlich and M. J. Taylor, \textit{Algebraic number theory},
Cambridge University Press, 1991.

[6] J. C. Lagarias, \textit{Meyer's concept of quasicrystal and quasiregular
sets}, Commun. Math. Phys. \textbf{179} (1996), 365-376.

[7] Y. Meyer, \textit{Nombres de Pisot, nombres de Salem et analyse
harmonique,} Lecture Notes in Mathematics, Springer-Verlag, \textbf{117}
(1970).

[8] R. V. Moody, \textit{Meyer sets and their duals,} The Mathematics of
long-range aperiodic order (Waterloo, ON, 1995), 403-441, NATO Adv. Sci.
Inst. Ser. C Math. Phys. Sci., \textbf{489}, Kluwer Acad. Publ., Dordrecht,
1997.

[9] C. Pisot, \textit{Quelques aspects de la th\'{e}orie des entiers alg\'{e}%
briques,} S\'{e}minaire de math\'{e}matiques sup\'{e}rieures, Universit\'{e}
de Montr\'{e}al, 1963.

[10] T. V\'{a}vra and F. Veneziano, \textit{Pisot units generators in a
number field}, J. Symb. Comp. \textbf{89 }(2018)\textbf{, }94-108.

[11] T. Za\"{\i}mi, \textit{On }$\varepsilon -$\textit{Pisot numbers}, New
York J. Math. \textbf{15 }(2009)\textbf{, }415-422.

[12] T. Za\"{\i}mi, \textit{Commentaires sur quelques r\'{e}sultats sur les
nombres de Pisot}, J. Th\'{e}or. Nombres Bordx. \textbf{22}, No. 2 (2010),
513-524.

[13] T. Za\"{\i}mi, \textit{On the distances between Pisot numbers
generating the same number field,} Bull. Malays. Math. Sci. Soc. \textbf{48}%
, No. 1 (2025), Paper No. 8, 15 p.

\bigskip

\bigskip

M. J. Bertin, Sorbonne Universit\'{e}, 4 Place Jussieu, 75252 Paris, Cedex
5, France

e-mail: marie-jose.bertin@imj-prg.fr\textit{\ }

\bigskip

T. Za\"{\i}mi, National Higher School of Mathematics, Sidi Abd-Allah,
Technology Hub, P. O. Box 75,\ \ Algiers 16093, Algeria

e-mail: toufik.zaimi@nhsm.edu.dz

\end{document}